\documentclass{amsart}[11pt]
\usepackage{amsmath, amscd,amsthm, amsfonts}
\usepackage[mathscr]{eucal}
\usepackage[all]{xy}
\usepackage[margin=2.3cm]{geometry}

\def\V{\mathbb V}
\def\R{\mathbb R}

\def\H{\mathbb H}
\def \h{\bf H}

\def\W{\mathbb W}
\def\C{\mathbb C}

\def\s{\mathbb S}

\def\F{\mathbb F}

\newtheorem{theorem}{Theorem}[section]
\newtheorem{lemma}[theorem]{Lemma}
\theoremstyle{definition}

\newtheorem{proposition}[theorem]{Proposition}
\theoremstyle{remark}

\numberwithin{equation}{section}
\theoremstyle{plain}

\newtheorem{corollary}[theorem]{Corollary}

\newcommand{\secref}[1]{section~\ref{#1}}
\newcommand{\thmref}[1]{Theorem~\ref{#1}}

\newcommand{\Propref}[1]{Proposition~\ref{#1}}

\newcommand{\eqnref}[1]{~{\textrm(\ref{#1})}}

\begin{document}
\title[Algebraic Characterization of isometries of the hyperbolic plane]{Algebraic Characterization
 of the Isometries of the  Complex and  Quaternionic Hyperbolic Plane}
\author[Wensheng Cao and Krishnendu Gongopadhyay]{Wensheng Cao  and  Krishnendu Gongopadhyay}
\address{School of Mathematics and Computational Science, Wuyi University,
Jiangmen, Guangdong 529020, P.R. China} \email{wenscao@yahoo.com.cn}
\thanks{During the research Cao has been supported by NSFS of China (10801107, 11071063) and
NSF of Guangdong Province (8452902001000043). }
\address{Indian Institute of Science Education and Research (IISER) Mohali, Phase 9, Sector 81, Mohali 160062, India}
\email{krishnendu@iisermohali.ac.in, krishnendug@gmail.com}
\date{Feb 23, 2011}
\begin{abstract}
Let $\h^2_{\F}$ denote two dimensional hyperbolic space over
$\F$, where $\F$  is either the complex numbers $\C$ or the quaternions $\H$. It
is of  interest to characterize algebraically 
the dynamical types of the isometries of $\h^2_{\F}$. For $\F=\C$, such a
characterization is known from the work of Giraud-Goldman. In this
paper,
 we offer an algebraic characterization of the isometries of $\h^2_{\H}$.
 Our result  restricts to the case  $\F=\C$ and provides another characterization of the isometries of $\h^2_{\C}$ which is different from the characterization due to Giraud-Goldman. 

Two elements in a group $G$ are said to be in the same \emph{$z$-class} if their centralizers are conjugate in $G$.
The $z$-classes provide a finite partition of the isometry group.
 In this paper we describe the centralizers  of the isometries of $\h^2_{\F}$,  and determine the $z$-classes.
\end{abstract}
\maketitle

\section{Introduction}
Let $\F$ denote the real numbers $\R$, the complex numbers $\C$ or the quaternions $\H$.
Let $\h^n_{\F}$ denote the $n$-dimensional hyperbolic space over $\F$. For $\F=\C$ and $\H$, 
 the linear groups which act as the isometries, are  denoted by $U(n,1)$ and $Sp(n,1)$ respectively. For $\F=\R$, we consider the group of orientation-preserving isometries of $\h^n_{\R}$,  and the linear group which acts as the orientation-preserving isometries is denoted  by $SO_o(n,1)$. The actual isometry groups  are the projectivizations of these linear groups.  For a group $G$,  recall that the projectivization of $G$, denoted by $PG$, is the quotient of $G$ by its center,  i.e. $PG=G/Z(G)$.

It is well known that the rank-one symmetric spaces of non-compact
type are the real, complex and quaternionic hyperbolic spaces,  and
the Cayley hyperbolic planes. Real hyperbolic geometry is extensively studied, complex hyperbolic geometry less so, but it is still a central subject of research (see \cite{bm},  \cite{gold},  \cite{parker}, \cite{sc}).  The quaternionic hyperbolic space is less well-understood.   Recently  an attempt has been made by Kim-Parker \cite{kp} to understand quaternionic hyperbolic geometry.

Classically,  one identifies the isometries of $\h^2_{\R}$ and $\h^2_{\C}$  with the invertible $2 \times 2$ matrices over the real and the complex numbers respectively. Then the dynamics of the isometries are classified in terms of the function   $\frac{trace^2}{det}$  ( see  \cite[Theorems 4.3.1 and 4.3.4]{bea},  \cite[Appendix-1]{gk}).  Algebraic characterization of the isometries of the real hyperbolic $4$-space can
 be obtained from  \cite{cpw} or  \cite{kido}.
  An alternative  characterization of the isometries of $\h^4_{\R}$ may also be obtained from   \cite[Theorem-1.1]{g} using the work of Ahlfors \cite{ahlfors} and Waterman \cite{waterman} in the quaternionic setting.  There are two remarkable Lie theoretic isomorphisms which identify
$PU(1,1)$ and $PSp(1,1)$ with the isometries of $\h^2_{\R}$ and
$\h^4_{\R}$, respectively. Therefore the isometries of $\h^1_{\C}$ and  $\h^1_{\H}$ can be characterized algebraically using results from low dimensional real hyperbolic geometry. 

It is an interesting question to ask for similar characterizations in higher dimensional complex and quaternionic hyperbolic geometries.  In the literature, Giraud \cite{giraud}  is known to be the first person who obtained  an algebraic characterization of the isometries of $\h^2_{\C}$.   Later,  Goldman \cite[Theorem 6.2.4]{gold}   developed a far more complete classification of the isometries in two dimensional complex hyperbolic geometry.  
   In this paper we offer an algebraic characterization of the isometries of $\h^2_{\H}$. This gives an effective algorithm to
detect the dynamical types of the isometries. Our result may be considered
as a counterpart of Goldman's theorem in the quaternionic hyperbolic
setting. As a corollary, we have obtained another characterization of the
 isometries of $\h^2_{\C}$. This is different from the characterizations due to Giraud-Goldman. 

 From now on, let $\F$  be either $\C$ or  $\H$.  When $\F=\H$, the underlying Hermitian space $\H^{2,1}$  is assumed to be a \emph{right} vector space over $\H$. Therefore, we consider only right eigenvalues of $g$ in $Sp(2,1)$.  

Let $g$ be an isometry of $\h^2_{\F}$. This lifts to a unitary transformation $\tilde g$ in the linear group of isometries. In the projective model of the hyperbolic space, the fixed points of $g$ correspond to the (right) eigenvectors of $\tilde{g}$. Since linear maps are easier to deal with, we shall mostly work with the lifts of the isometries.  Therefore $g$ will be identified with $\tilde g$, and forgetting `tilde'  we shall denote them by the same symbol. 

When $\F=\H$, the eigenvalues of an isometry occur in similarity classes, i.e. if $\lambda$ is an eigenvalue of $g$, then $v \lambda v^{-1}$ is also an eigenvalue,  where $v \in \H-\{0\}$.  Therefore, an eigenvalue of $g$ will be understood as a similarity class of eigenvalues.  We say \emph{two eigenvalues of $g$ are distinct} if the corresponding similarity classes of the eigenvalues are disjoint.  Two eigenvalues are \emph{equal} if they belong to the same similarity class.  Each similarity class of eigenvalues contains a unique pair of complex conjugate numbers.
 We further adopt the convention of choosing the eigenvalue $re^{i \theta}, 0 \leq \theta \leq \pi$ from the similarity class,  and we identify the similarity class with this complex eigenvalue.

 By the Brouwer's fixed-point theorem, every isometry of $\h^2_{\F}$ has a fixed point on
$\h^2_{\F} \cup \partial \h^2_{\F}$. An isometry $g$ is \emph{elliptic} if it has a fixed point on $\h^2_{\F}$;  \emph{parabolic} if it is non-elliptic and has a \emph{unique}  fixed point on $\partial \h^2_{\F}$;   \emph{hyperbolic} if it is non-elliptic and has exactly \emph{two } fixed points on  $\partial \h^2_{\F}$. We refine this classification further as follows. 

\begin{itemize}
\item[(e)]{  Let $g$ be elliptic:  If $g$ has mutually distinct eigenvalues, then $g$ is  called a  \emph{regular elliptic}. If $g$ has two equal eigenvalues,  it is called  a \emph{complex elliptic}, (or a \emph{complex reflection}). If all the eigenvalues of $g$ are equal then we call it a \emph{simple elliptic}. 
 Note that simple elliptics occur only if $\F=\H$. In the complex case, they belong to the center of the group and hence act as the identity.}

 \item[(h)]{  Suppose $g$ is hyperbolic:  It follows from the conjugacy classification that it has a complex eigenvalue outside the unit disc, and one eigenvalue inside the unit disc. The other eigenvalue lies on the unit circle. The isometry $g$ is called a \emph{regular hyperbolic} if it has a non-real eigenvalue of norm different from $1$.  If all the eigenvalues of $g$ are real numbers, then it is called a  \emph{strictly hyperbolic}.  If $g$ has two and only two real eigenvalues, then it is called  a \emph{screw hyperbolic}.   Note that screw hyperbolics occur only if $\F=\H$. When $\F=\C$,  dynamically there is no difference between the regular hyperbolics and the screw hyperbolics.  It can be seen from the conjugacy classification in $U(2,1)$ that a regular hyperbolic can be obtained from a screw hyperbolic by multiplying a central element of the form  $\lambda I$, where $|\lambda|=1$, and vice-versa.  Hence we identify the two classes in $U(2,1)$, and following Goldman,  call them \emph{loxodromic}.}

 \item[(p)] { Suppose $g$ is parabolic:   Let $g$ be unipotent, i.e.  all eigenvalues of $g$ are $1$.  If the minimal polynomial of $g$  is $(x-1)^2$, then it is called a \emph{vertical   translation}. It  is a \emph{non-vertical   translation} if the minimal polynomial is $(x-1)^3$. Suppose $g$ is a  non-unipotent, i.e. it has a non-real
  eigenvalue. Suppose the multiplicity of the non-real eigenvalue is $3$:  then $g$ is an  \emph{ellipto-translation} or \emph{ ellipto-parabolic} according as the minimal polynomial of $g$ has degree $2$ or $3$.  When $\F=\C$, these classes do not occur, as they act on
$\h^2_{\C}$ as the vertical or non-vertical   translations. If $g$ has two distinct eigenvalues,  it is called a \emph{screw parabolic}.  }
\end{itemize}
Given any isometry $g$ of $\h^2_{\F}$, it belongs to one of the above classes.

We use  the embedding of $\H$ into $2
\times 2$ complex matrices $M_2(\C)$.
 This gives an embedding $A \mapsto A_{\C}$ of $Sp(2,1)$ into $ GL(6, \C)$. Using this embedding we obtain the following characterization of the isometries.
 \begin{theorem}\label{mainth}
Let $A$ be an element in $Sp(2,1)$. Let $A_{\C}$ be the
corresponding element in $GL(6, \C)$. The characteristic polynomials
of $A_{\C}$ is of the form
\begin{equation*}\chi_A(x)=x^6-ax^5 +bx^4-cx^3+bx^2-ax+1,\end{equation*}
where $a$, $b$, $c$ are real numbers.
Define $$G=27(a-c)+9ab-2a^3,$$
$$H=3(b-3)-a^2,$$
$$\Delta=G^2 +4H^3.$$
Then we have the following.

1.     $A$ acts as a regular hyperbolic if and only if $\Delta >0$.

2. $A$ acts as a regular elliptic if and only if $\Delta<0$. 

3.  $A$ acts as  either strictly hyperbolic, screw hyperbolic,  complex elliptic,  or  screw parabolic if and only if $\Delta=0$ and $G \neq 0$.
Moreover, we have the following.   
\begin{itemize}
\item[]{ Let $\alpha_n=trace(A_{\C}^n)$. }
\item[(i)]{ $A$ acts as either strictly hyperbolic or screw hyperbolic if and only if the sequence $\{\alpha_n\}$ is divergent.  Further $A$ acts as a strictly hyperbolic if and only if $(2a+c)^2=4(1+b)^2$.   }

\item[(ii)]{ $A$ acts as either complex elliptic or  screw parabolic if and only if the sequence $\{\alpha_n\}$ is bounded.  Further $A$ acts as a screw parabolic if and only if the degree of the minimal polynomial of $A_{\C}$ is $3$.    }
\end{itemize}

4. $A$ acts as either simple elliptic,  ellipto-translation or ellipto-parabolic if and only if $\Delta=0$, $G=0$ and $|a|<6$, $|b|<15$, $|c|<20$. 
Further, $A$ acts as a simple elliptic, resp. ellipto-translation, resp. ellipto-parabolic if and only if the degree of the minimal polynomial
of  $A_{\C}$, is $1$,  resp. $2$, resp. $3$. 

5. $A$ acts as a unipotent if and only if  $A \neq \pm I$ and $|a|=6, \ |b|=15,  \ |c|=20$. Further,  $A$ acts as a vertical, or a non-vertical translation according as the  degree of the minimal polynomial of $A_{\C}$  is $2$, or $3$. 
\end{theorem}
For each $n$,  $\alpha_n$  is the sum of $n$-th powers of  eigenvalues of  $A_{\C}$.  The well-known Newton's identities  (see  \cite[Theorem-1.3.19]{row}, \cite{mcd})    express the $\alpha_n$'s in terms of the co-efficients $a, \ b, \ c$ of the characteristic polynomial.  
One advantage of our method is that it applies also to $U(2,1)$. In this case the above embedding restricts to the embedding $A \mapsto A_{\R}$ of $U(2,1)$ into $GL(6, \R)$. This provides the following characterization of the isometries of $\h^2_{\C}$. 
\begin{corollary}\label{cor1}
 Let $A$ be an element in $U(2,1)$ Let $A_{\R}$ be the corresponding element in $GL(6, \R)$. The characteristic polynomials of $A_{\R}$ is of the form
\begin{equation*}\chi_A(x)=x^6-ax^5 +bx^4-cx^3+bx^2-ax+1.\end{equation*}
Define $$G=27(a-c)+9ab-2a^3,$$
$$H=3(b-3)-a^2,$$
$$\Delta=G^2 +4H^3.$$
$$\alpha_n=trace(A_{\R}^n)$$
Then we have the following.

1.   $A$ acts as a regular elliptic if and only if $\Delta<0$.

2.  $A$ acts as a loxodromic if and only if one of the following conditions hold. 
\begin{itemize}
\item[ (i)]{ $\Delta>0$, }  
\item[(ii)]{ $\Delta=0$, $G \neq 0$,  the sequence $\{\alpha_n\}$ is divergent and $(2a+c)^2 \neq 4(1+b)^2$. }
\end{itemize}
3. $A$ acts as a strictly hyperbolic if and only if $\Delta=0$, $G \neq 0$,   the sequence $\{\alpha_n\}$ is divergent, and  $(2a+c)^2 = 4(1+b)^2$.

4.   $A$ acts as either complex elliptic or,  screw parabolic  if and only if   $\Delta=0$, $G \neq 0$ and  the sequence $\{\alpha_n\}$ is bounded. Further, $A$ acts as a screw parabolic if and only if the degree of the minimal polynomial of $A$, or $A_{\R}$, is $3$.

5. $A$ acts as a unipotent if and only if $\Delta=0$, $G=0$, 
$|a| \leq 6$, $|b|\leq 15$, $|c| \leq 20$. It is a vertical, or a non-vertical   translation according as the degree of the minimal polynomial is $2$ or $3$.

6.  $A$ acts as the identity if and only if $A=\lambda I$, where $|\lambda|=1$. 
\end{corollary}

  Two elements $x$, $y$ in a group $G$ are said to
be in the same $z$-class if their centralizers are conjugate in $G$.
This notion, and terminology, are due to Kulkarni \cite{rkrjm}. The
$z$-classes in a group give a partition of $G$, and they refine the
partition of $G$ into conjugacy classes. In the case of compact Lie
groups, and many other cases, there are only {\it {finitely many}}
$z$-classes. The notion is defined in terms of the group-structure
alone. Kulkarni \cite{rkrjm} has proposed that this notion
may be used to make precise the intuitive idea of ``dynamical
types'' in any ``geometry'', whose automorphism group contains a
copy of $G$.
 If the $z$-classes are finitely many in number, then they
provide a  finite classification of the group of transformations. Thus the partition of the group into its $z$-classes provide another algebraic
 characterization of the isometries, and this is based on the internal structure of the isometry group alone.

Let $PO(n,1)$ denote the full group of isometries of the real hyperbolic space $\h^n_{\R}$. Using the linear model of the hyperbolic space, the 
$z$-classes in $PO(n,1)$  are classified in \cite{gk}. In
particular,  this classification is valid when $n=2$.   Classification of centralizers of the orientation-preserving isometries of $\h^2_{\R}$ can also be obtained from \cite[Lemma-1.10]{scott}. 
In this paper, we classify the centralizers, and the $z$-classes
in $U(2,1)$ and $Sp(2,1)$.  In fact, we compute their number. 
\begin{theorem}\label{zcl}
(1)  In the group $Sp(2,1)$, there are eleven $z$-classes
of elliptic elements, four $z$-classes of hyperbolic elements, two
$z$-classes of unipotent elements and six $z$-classes of parabolic non-unipotent elements. In this count we consider the
identity map as an elliptic element.

(2)  In the group $U(2,1)$, there are four
$z$-classes of elliptic elements, a unique $z$-class of hyperbolic
elements, two $z$-classes of unipotent elements and a unique $z$-class of parabolic non-unipotent elements. In this count we
consider the identity map as an elliptic element.
\end{theorem}

 To count the $z$-classes in $PO(n,1)$, a combinatorial formula  is available  in  \cite[Theorem-1.3]{gk}.  Such formula is not used in the proof of the above theorem. Instead, the proof of the theorem relies on case-by-case study of the centralizers of the isometries.  A related problem is to classify the commuting isometries. The commuting elements in $SU(2,1)$ are classified by Basmajian-Miner \cite{bm}. Recently, the authors \cite{cg2} have generalized this classification in arbitrary dimensions.   Along the way,  the $z$-classes in $U(n,1)$ are also classified. Using similar methods,  the $z$-classes in $Sp(n,1)$ may be classified  as well.    However, counterparts of \cite[Theorem-1.3]{gk} for the groups $U(n,1)$ and $Sp(n,1)$  are still missing in the literature. 

The structure of the remainder of this paper is as follows. In
 \secref{preli}, we note some preliminary results. The conjugacy classes are crucial in the proof of the above theorems, and they are classified in \secref{cc}. We prove our main theorem, viz. \thmref{mainth},  in \secref{pfmainth}. We classify the centralizers and prove \thmref{zcl} in \secref{zcc}.

\section{Preliminaries}\label{preli}
\subsection{The Hyperbolic Space} Assume $\F$ to be either  $\C$ or $\H$. When considering a vector space over $\H$, it is assumed that the scalar multiplication acts from the \emph{right}, hence the respective subspaces are \emph{right subspaces}. 

 Let $\V=\V^{n,1}$ be the right vector space of dimension $n+1$  over $\F$ equipped with a  non-degenerate Hermitian form $\langle, \rangle$ of signature $(n,1)$. With respect to a suitable basis, the Hermitian form $\langle, \rangle$ is given by 
$$
\langle{\bf z},\,{\bf w}\rangle={\bf z}^*J_1{\bf w}=
-\overline{z_0}w_0+\overline{z_{1}}w_{1}+\overline{z_{2}}w_{2}+...+\overline{z_n}w_n,
$$
where  $\cdot^*$   denotes the Hermitian transpose,  ${\bf z}$ denotes the column vector in $\V$ with entries $z_0,\ z_1,\  z_2,.....,z_n$,  and
\hbox{$J_1=diag(-1,1, 1,...,1)$}.   An \hbox{\emph{ isometry}} $g$ of $\V$  is  a (right) linear
bijection such that  for all ${\bf z}$ and
${\bf w}$ in $\V$, $\langle g({\bf z}),\,g({\bf
w})\rangle=\langle{\bf z},\,{\bf w}\rangle$. The  isometry group is denoted by $U(n,1;\F)$.
The group of isometries of the Hermitian form
$$Q(z,w)=\overline{z_{1}}w_{1}+\overline{z_{2}}w_{2}+...+\overline{z_{n}}w_{n}$$
is denoted by $U(n;\F)$. Abbreviating the symbols,  we denote $U(n,\C)=U(n)$, $U(n,
\H)=Sp(n)$, $U(n,1;\C)=U(n,1)$ and $U(n,1;\H)=Sp(n,1)$.

Our interest
in this paper is the case when $n=2$.
Following Section 2 of  \cite{chen}, let
\begin{eqnarray*}
\V_0 & = & \Bigl\{{\bf z} \in  \V^{2,1}-\{0\}:
\langle{\bf z},\,{\bf z}\rangle=0\Bigr\} \\
\V_{-} &  = & \Bigl\{{\bf z} \in \V^{2,1}:\langle{\bf z},\,{\bf
z}\rangle<0\Bigr\}.
\end{eqnarray*}
It is obvious that $\V_0$ and $\V_{-}$ are invariant under
$U(2,1;\F)$. Let $\V^s=\V_{-}\cup  \V_0$. Let
$P:\V^s\to P(\V^s)\subset \V^{2}$ be the right  projection map
defined by
$$
P(z_0, z_1, z_{2})^t=(z_1z_{0}^{-1},  z_2z_{0}^{-1})^t,
$$
where $\cdot^t$ denotes the  transpose. We define $\h_{\F}^2=P(\V_-)$. The boundary of the hyperbolic space is $\partial \h_{\F}^2=P(\V_0)$.

It is convenient to introduce the Cayley transform
$$C=\left(
      \begin{array}{ccc}
        \sqrt{2}/2 & -\sqrt{2}/2& 0\\
       \sqrt{2}/2&\sqrt{2}/2 & 0 \\
        0 & 0 & 1\\
      \end{array}
    \right).$$

The Cayley transformation $C$
maps $\h_{\F}^2$ and its boundary $\partial \h_{\F}^2$ to the Siegel
domain $\Sigma_{\F}^2$ and its boundary $\partial \Sigma_{\F}^2$,
respectively.  The Hermitian form transforms to
$$
\langle{\bf z},\,{\bf w}\rangle=
-(\overline{z_0}w_1+\overline{z_{1}}w_{0})+\overline{z_{2}}w_{2}.
$$
The isometry group of the above form is $\hat U(2,1;\F)=CU(2,1;\F)C^{-1}$,

\begin{proposition}
(1) If $A=\left(
            \begin{array}{ccc}
              a & b & c \\
              d & e & f \\
              g & h & l \\
            \end{array}
          \right)\in U(2,1;\F)$,  then $A^{-1}=\left(
            \begin{array}{ccc}
              \bar{a} & -\bar{d} & -\bar{g} \\
              -\bar{b} & \bar{e} & \bar{h} \\
              -\bar{c} & \bar{f} & \bar{l} \\
            \end{array}
          \right).$

(2) If $A=\left(
            \begin{array}{ccc}
              a & b & c \\
              d & e & f \\
              g & h & l \\
            \end{array}
          \right)\in \hat U(2,1;\F)$,  then $A^{-1}=\left(
            \begin{array}{ccc}
              \bar{e} & \bar{b} & -\bar{h} \\
              \bar{d} & \bar{a} & -\bar{g} \\
              -\bar{f} & -\bar{c} & \bar{l} \\
            \end{array}
          \right).$
\end{proposition}

Let  $o=C(-f_1), \ \infty=C(f_1)\in \partial \Sigma$, where
$f_1=(1,0)\in \partial \h^2_{\F}$.  Let
$$G_0=\{g\in \hat U(2,1;\F): g(o)=o\}, \ G_{\infty}=\{g\in \hat U(2,1;\F):
g(\infty)=\infty\}, \ G_{0,\infty}=G_0\cap G_{\infty}.$$ Then we have
the following three propositions, cf.   \cite[Lemma 3.3.1]{chen}.

\begin{proposition}
If $A\in G_{\infty}$, then $A$ is of the form
 \begin{equation}\label{gw}\left(
  \begin{array}{ccc}
    a & 0 & 0 \\
    d & e & f \\
    g & 0 & l \\
  \end{array}
\right), \ \mbox{where}\ |l|=1,\
\bar{a}e=1,\;\Re(\bar{a}d)=\frac{1}{2}\left|g\right|^2,\;f=e\bar{g}l.\end{equation}
\end{proposition}

\begin{proposition}
If $A\in G_{0}$, then $A$ is of the form
 \begin{equation}\label{g0}\left(
  \begin{array}{ccc}
    a & b & c \\
    0 & e & 0 \\
    0 & h & l \\
  \end{array}
\right), \ \mbox{where}\ |l|=1,\
\bar{a}e=1,\;\Re(\bar{e}b)=\frac{1}{2}\left|h\right|^2,\;c=a\bar{h}l.\end{equation}
\end{proposition}

\begin{proposition}
If $A\in G_{0,\infty}$, then $A$ is of the form
 \begin{equation}\label{g0w}\left(
  \begin{array}{ccc}
    a & 0 & 0 \\
    0 & e & 0 \\
    0 & 0 & l \\
  \end{array}
\right), \ \mbox{where}\   |l|=1,\ \bar{a}e=1.\end{equation}
\end{proposition}

\begin{lemma}\label{type1}
Let $A\in \hat U(2,1;\F)$ be of the form (\ref{parae}). Then\\
(i) $A$ acts as an  elliptic isometry if and only if
\begin{equation*}\label{rke}
rank\left(\begin{array}{cc}
f & d\\
 e^{{\bf i}  \phi}-e^{{\bf i}  \theta} & g \\
\end{array}\right)= rank \left(\begin{array}{c}
 f\\ e^{{\bf i}  \phi}-e^{{\bf i}  \theta}
\end{array}\right);
\end{equation*}\\
(ii)  $A$ acts as a parabolic isometry if and only if
\begin{equation*}\label{tobepara}
rank\left(\begin{array}{cc}
f & d\\
 e^{{\bf i}  \phi}-e^{{\bf i}  \theta} & g \\
\end{array}\right)\neq rank \left(\begin{array}{c}
 f\\ e^{{\bf i}  \phi}-e^{{\bf i}  \theta}
\end{array}\right).
\end{equation*}
\end{lemma}

\begin{proof}
 Let $g$ be the isometry corresponding to $A$ in the Siegel domain model.
The isometry  $g$ acts on $\overline{\Sigma_{\F}^2}$  as follows.
\begin{equation*}
g\left(
\begin{array}{c}
               \eta_1 \\
               \eta_2 \\
             \end{array}
           \right)= e^{-{\bf i}\theta}\left(
                    \begin{array}{cc}
                      e^{{\bf i}  \theta} & f \\
                      0 & e^{{\bf i}  \phi} \\
                    \end{array}
                  \right)\left(
\begin{array}{c}
               \eta_1 \\
               \eta_2 \\
             \end{array}
           \right)
+e^{-{\bf i}\theta}\left(
             \begin{array}{c}
               d \\
               g \\
             \end{array}
           \right).
\end{equation*}
 Thus $g$ has a fixed point $\eta=(\eta_1, \eta_2)^t$ if and only if the following set of equations has a solution.
\begin{equation}\label{rkeqe1} \left(
                    \begin{array}{cc}
                      0  & f \\
                      0 & e^{{\bf i}\phi}-e^{{\bf i}\theta} \\
                    \end{array}
                  \right)\left(
\begin{array}{c}
               \eta_1 \\
               \eta_2 \\
             \end{array}
           \right)=-\left(
             \begin{array}{c}
               d \\
               g \\
             \end{array}
           \right).
 \end{equation}
It follows from basic linear algebra that  \eqnref{rkeqe1} has a solution if and only if the condition in $(i)$ holds.
Clearly $g$ is elliptic if and only if \eqnref{rkeqe1} has a solution. Otherwise, it is parabolic.

  The proof is now complete.
\end{proof}

\begin{lemma}
The group $Sp(2,1)$ can be embedded in the group $GL(6, \C)$.
\end{lemma}
The proof is similar to that of Proposition-2.4 in  \cite[p.160]{g}, also see
 \cite[section-2]{lee}  or  \cite[section-2]{zhang}.   We sketch the proof here.
\begin{proof}
Write $\H=\C\oplus{\bf j}  \C$. For $A\in Sp(2,1)$, express $A=A_1+{\bf j}A_2$,
{where} $ A_1, A_2\in M_3(\C)$. This gives an embedding $A \mapsto A_{\C}$ of $Sp(2,1)$ into $GL(6, \C)$, where
\begin{equation}\label{crep} A_{\C}= \left(
                          \begin{array}{cc}
                            A_1 &  -\overline{A_2}\\
                          {A_2}  & \overline{A_1} \\
                          \end{array}
                        \right).
\end{equation}
\end{proof}
Similarly, for $A \in SU(2,1)$, one can express $A=A_1+A_2{\bf i}, \
\mbox{where}\ A_1, A_2\in M_3(\R)$ and we have
\begin{equation}\label{rrep} A_{\R}= \left(
                          \begin{array}{cc}
                            A_1 & -{A_2} \\
                           {A_2} & {A_1} \\
                          \end{array}
                        \right).
\end{equation}

\subsection{Preliminaries on roots of polynomials}\label{prel}
Here we note down a few facts about the nature of solutions of real cubic equations. 
For details see \cite{n1, n2}.

\subsubsection{Resultant of two polynomials}
Let $R(f,g)$ denote the resultant of two polynomials $f(x)$ ad $g(x)$ over the reals.  The \emph{resultant} of a polynomial $f(x)$ is,
 up to a scalar factor, the resultant $R(f, f')$. Recall that the polynomial $f(x)$ has a multiple root if and only if the resultant is zero.
  The polynomial has a triple root if and only if the resultant $R(f, f'')$ is zero, and so on.

\subsubsection{Nature of solutions of a cubic equation}\label{nature}
Let $f(x)=ax^3 + 3bx^2+3cx+d$ be a polynomial over the reals.
Let $G=a^2d-3abc+2b^3$, $H=ac-b^2$, and $\Delta=G^2+4H^3$. Then the type of the roots of the equation $f(x)=0$ can be detected as follows.

(i) If $\Delta>0$, then only root of the equation is real, and the other two are complex conjugates.

(ii) If $\Delta<0$, then all the roots of the equation are real and distinct.

(iii) If $\Delta=0$, then the roots are real, and at least two of them are equal. The equation has a triple real
root if and only if the resultant $R(f,  f'')$ vanishes.

\section{The Conjugacy Classification}\label{cc}
\begin{theorem}\label{conjugation}
 (i) Suppose $A\in \hat U(2,1;\F)$ acts as a hyperbolic element. Then $A$ is conjugate to an isometry  of the form
\begin{equation}\label{loxoe} L=L(\beta, \theta)=\begin{pmatrix} re^{{\bf i} \beta} & 0 & 0 \\
0 & r^{-1}e^{{\bf i}  \beta} & 0 \\0 & 0 & e^{{\bf i}  \theta}
\end{pmatrix}, \;r>0,\;r\neq 1,\end{equation}
when $\F=\H$, $0 \leq \beta, \theta \leq \pi$. For $\F=\C$, $-\pi \leq \beta, \  \theta \leq \pi$.

(ii)  Suppose $A\in U(2,1;\F)$ acts as an elliptic element. Then  $A$ is conjugate to an isometry of the form
\begin{equation}\label{ellie}E=E(\theta, \phi, \psi)=\begin{pmatrix} e^{{\bf i}  \theta} & 0 & 0 \\ 0 & e^{ {\bf i}  \phi} & 0 \\0 & 0 & e^{{\bf i}  \psi}
\end{pmatrix},\end{equation}
when $\F=\H$, $0\leq \theta, \phi, \psi \leq \pi$, and when $\F=\C$, $-\pi \leq \theta, \phi, \psi \leq \pi$.

(iii) Suppose $A\in \hat U(2,1;\F)$ acts as a  parabolic element. Then $A$ is conjugate to an isometry of the form
\begin{equation}\label{parae}P=\left(
  \begin{array}{ccc}
    e^{{\bf i}\theta} & 0 & 0 \\
    d & e^{{\bf i}\theta} & f \\
    g & 0 & e^{{\bf i}\phi} \\
  \end{array}
\right),\end{equation}   where $d \in \C-\{0\}$, $f, g\in \C, \ \Re(e^{-{\bf i}\theta}d)=\frac{1}{2}|g|^2,  f=e^{{\bf
i}(\theta+\phi)}\bar{g}$, $0 \leq \theta, \ \phi \leq \pi$,  $(\hbox{resp. } -\pi \leq \theta, \ \phi \leq \pi)$ for $\F=\H$, $(\hbox{resp. }\C)$,
and
\begin{equation*}\label{tobepara1}
rank\left(\begin{array}{cc}
f & d\\
 e^{{\bf i}  \phi}-e^{{\bf i}  \theta} & g \\
\end{array}\right)\neq rank \left(\begin{array}{c}
 f\\ e^{{\bf i}  \phi}-e^{{\bf i}  \theta}
\end{array}\right).
\end{equation*}
\end{theorem}

\medskip  We shall prove the theorem over $\H$ following \cite[section-3]{chen}.  The proof over $\C$ is completely  analogous.   
\begin{proof}
If
$A\in \hat U(2,1;\H)$ is loxodromic, then $A$ has two fixed points
on $\partial \Sigma_{\H}^2$.  By conjugation, if necessary, we may assume that $A$ is of the form (\ref{g0w}). Let $u,v$  be quaternions of unit modulus such
that $l=ve^{{\bf i}\theta}v^{-1}$ and $a=ure^{{\bf i}\beta}u^{-1}$,
where $r>0, r\neq 1, 0\leq \theta, \beta \leq \pi$.  Let $U=diag(u,u,v)\in \hat U(2,1;\H)$. Then we have
\begin{equation*}
UAU^{-1}=diag(re^{{\bf i}\beta},r^{-1}e^{{\bf i}\beta},e^{{\bf
i}\theta}),
\end{equation*}
 This completes the proof of $(i)$.

By \cite[Proposition 3.2.1]{chen} and \cite[Corollary 6.2]{zhang},
if $A\in Sp(2,1)$ is elliptic, then $A$ is conjugate to
$$E=E(\theta, \phi, \psi)=diag(e^{{\bf i}  \theta}, e^{ {\bf i}  \phi}, e^{{\bf i}
\psi}),\  0\leq \theta, \phi, \psi \leq \pi.$$
This completes the proof of $(ii)$.

For the parabolic isometries, we have the following cases.

{\bf Case 1:}  All the right eigenvalues of $A$ are similar.

Suppose $A$ has no real right eigenvalue. Following
Chen-Greenberg \cite[Lemma 3.4.2]{chen}, $A$ is conjugate to
\begin{equation*}\label{gw1}\left(
  \begin{array}{ccc}
     e^{{\bf i}\theta} & 0 & 0 \\
    d &  e^{{\bf i}\theta} & f \\
    g & 0 &  e^{{\bf i}\theta} \\
  \end{array}
\right)\in \hat U(2,1;\C),\end{equation*} where $\Re( e^{-{\bf
i}\theta}d)=\frac{1}{2}\left|g\right|^2,\;f= e^{2{\bf
i}\theta}\bar{g},\  0< \theta<
                \pi.$

                If   $A$ is parabolic with real right eigenvalues
                then it is conjugate to $$T=\left(
                  \begin{array}{ccc}
                    \lambda & 0 & 0 \\
                     d& \lambda & f \\
                   \bar{f}& 0 & \lambda \\
                  \end{array}
                \right), \ \ \mbox{where}\ \  \lambda=\pm 1, \Re(\lambda d)=\frac{1}{2}|f|^2, f,d\in\H.$$

Let $u\in\H$ be such that $|u|=1$ and $udu^{-1}\in \C$.
 Let
$uf=k_1+k_2{\bf j}$ with $k_1,k_2\in \C$ and $v=\left\{
  \begin{array}{ll}
    \frac{k_1}{|f|}+\frac{k_2}{|f|}{\bf j}, & \hbox{$f\neq 0$;} \\
    1, & \hbox{$f=0$.}
  \end{array}
\right.$

 Then
\begin{equation*}
BTB^{-1}=\left(\begin{array}{ccc}
                    \lambda & 0 & 0 \\
                    udu^{-1}& \lambda & ufv^{-1} \\
                    v\bar{f}u^{-1}& 0 & \lambda \\
                  \end{array}
                \right)\in \hat U(2,1;\C),
\end{equation*}
where $B=diag(u,u, v)\in \hat U(2,1;\H)$.

{\bf Case 2:} The right eigenvalues of $A$ fall into two
similarity classes. Then $A$ is conjugate to
\begin{equation*}\label{gw2}T=\left(
  \begin{array}{ccc}
     e^{{\bf i}\theta} & 0 & 0 \\
    d &  e^{{\bf i}\theta} & f \\
    g & 0 &  e^{{\bf i}\phi} \\
  \end{array}
\right),\end{equation*} where $\Re( e^{-{\bf
i}\theta}d)=\frac{1}{2}\left|g\right|^2,\;f= e^{{\bf
i}\theta}\bar{g}e^{{\bf i}\phi},\ d, g, f \in \H, 0\leq
\theta,\phi\leq
                \pi$ and $\theta\neq \phi$.

 Let $f=f_1+f_2{\bf
j},$ where $f_1,f_2\in \C$ and $c=\frac{f_2}{e^{-{\bf
i}\phi}-e^{{\bf i}\theta}}{\bf j}$, $q=\frac{|f_2|^2}{2|e^{-{\bf
i}\phi}-e^{{\bf i}\theta}|^2}$.  Then
\begin{equation*}T_1=B_1TB_1^{-1}=\left(
  \begin{array}{ccc}
     e^{{\bf i}\theta} & 0 & 0 \\
    s &  e^{{\bf i}\theta} & f_1 \\
    e^{{\bf i}\phi} \bar{f_1}e^{{\bf i}\theta} & 0 &  e^{{\bf i}\phi} \\
  \end{array}
\right),\ \ \mbox{where}\ \ B_1=\left(
  \begin{array}{ccc}
    1 & 0 & 0 \\
    q & 1 & \bar{c} \\
    c & 0 & 1\\
  \end{array}
\right)\end{equation*} and  $s$ may be not a complex number. Thus up to conjugacy, every element is of the form (\ref{parae}). Since $e^{i \theta} \neq e^{i \phi}$,
an element of the form  (\ref{parae}) is conjugate to
\begin{equation}\label{screwp}P=\left(
  \begin{array}{ccc}
    e^{{\bf i}\theta} & 0 & 0 \\
    d & e^{{\bf i}\theta} & 0 \\
    0& 0 & e^{{\bf i}\phi} \\
  \end{array}
\right).\end{equation} The conjugation is obtained by the map  $T=\left(
  \begin{array}{ccc}
     1 & 0 & 0 \\
    t &  1 & \bar{c} \\
    c & 0 &  1 \\
  \end{array}
\right)$,  where $\bar{c}=\frac{f}{e^{{\bf i}\theta}-e^{{\bf
i}\phi}}, t=\frac{1}{2}|c|^2$. The proof now follows from Lemma \ref{type1}.
\end{proof}

\section{Proof of Theorem \ref{mainth}}\label{pfmainth}
The following proposition follows from the conjugacy classification.
\begin{proposition}\label{keypo}
  Let  $A\in Sp(2,1)$, resp. $SU(2,1)$. The characteristic polynomial of  $A_{\C}$, resp. $A_{\R}$,   is of the form
$$\chi_A(x) = x^6-a_5x^5+a_4x^4-a_3x^3+a_2x^2-a_1x+a_0,$$
where $a_i\in \R$  and $a_1=a_5$, $a_2=a_4$, $a_0=1$.
\end{proposition}
\begin{proof}
Note that the characteristic polynomial is invariant under conjugation. Hence it is sufficient to consider $A_{\C}$ (resp. $A_{\R}$)
up to conjugacy. From the conjugacy classification, cf. \thmref{conjugation}, we observe the following.

\subsubsection*{ Hyperbolic isometries.}\label{hyp1}  Suppose $A$ is
conjugate to a matrix of the form (\ref{loxoe}) and the
characteristic polynomial of $A_{\C}$ is
\begin{equation*}
\chi_A(x)=(x^2-2\cos\theta\  x+1)(x^2-2r\cos\beta\
x+r^2)(x^2-\frac{2}{r}\cos\beta\ x+\frac{1}{r^2}).
\end{equation*}
Hence
\begin{equation}\label{loa1}
a_0=1,\  a_1=a_5=2(r+\frac{1}{r})\cos\beta+2\cos\theta,
\end{equation}
\begin{equation}\label{loa2}
a_2=a_4=4(r+\frac{1}{r})\cos\theta cos\beta+4\cos^2\beta
+r^2+\frac{1}{r^2}+1,
\end{equation}
\begin{equation}\label{loa3}
a_3=4(r+\frac{1}{r})\cos\beta+
2(r^2+\frac{1}{r^2}+4\cos^2\beta)\cos\theta.
\end{equation}

\subsubsection*{ Elliptic isometries}\label{ell1}  Suppose $A$ is
conjugate to an element of the form (\ref{ellie}) and the
characteristic polynomial of $A_{\C}$ is
\begin{equation*}
\chi_A(x)=(x^2-2\cos\theta\ x+1)(x^2-2\cos\phi\ x+1)(x^2-2\cos\psi\ x
+1)
\end{equation*}
 Hence
\begin{equation*}
a_0=1,\ a_1=a_5=2(\cos\theta+\cos\phi+ \cos\psi),
\end{equation*}
\begin{equation*}
a_2=a_4=3+4(\cos\theta\cos\phi +\cos\phi \cos\psi+\cos\psi
\cos\theta),
\end{equation*}
\begin{equation*}
a_3=4(\cos\theta+\cos\phi+ \cos\psi)+8\cos \theta \cos\phi \cos\psi.
\end{equation*}

\subsubsection*{ Parabolic isometries}\label{par1}  Suppose  $A$ is
conjugate to a matrix of the form (\ref{parae}) and the
characteristic polynomial of $A_{\C}$ is
\begin{eqnarray*}
\chi_A(x)&=&(x^2-2\cos\theta\ x+1)^2(x^2-2\cos\phi\ x+1).
\end{eqnarray*}
 Hence
\begin{equation*}
a_0=1,\ a_1=a_5=2(2\cos\theta+\cos\phi),
\end{equation*}
\begin{equation*}
a_2=a_4=3+4\cos^2\theta +8\cos\theta \cos \phi,
\end{equation*}
\begin{equation*}
a_3=4(\cos\phi+2\cos\theta)+8\cos^2 \theta \cos\phi.
\end{equation*}
This completes the proof of \Propref{keypo}.
\end{proof}

\subsection*{The Proof of \thmref{mainth}}  Observe that
$$\chi_A(x)=x^6-ax^5 +bx^4-cx^3+bx^2-ax+1$$ is a \emph{self-dual}
polynomial over $\C$, i.e. if $\alpha$ is a root of $\chi_A(x)$ then
$\alpha^{-1}$ is also a root. Hence we can write $\chi_A(x)=x^3
g_A(x)$, where
$$g_A(x)=(x^3 + x^{-3})-a(x^2 + x^{-2})+b(x+x^{-1})-c.$$
Further observe that $\chi_A(x)$ is, in fact, a polynomial over the reals. Hence if $\lambda$ in $\C$ is a root of $\chi_A(x)$, then $\bar \lambda$ is also a root.

In the expression of $g_A(x)$, expanding the terms in brackets,  we have
\begin{eqnarray*}
g_A(x)&=& (x+x^{-1})^3-3(x+x^{-1})-a[(x+x^{-1})^2-2]+b(x+x^{-1})-c.
\end{eqnarray*}
Let $t=x+x^{-1}$. Then $g_A(x)$ is a cubic polynomial in $t$, and we denote it by $g_A(t)$, i.e.
\begin{equation}\label{1}
g_A(t)= t^3-at^2+(b-3)t-(c-2a).
\end{equation}
Since $\chi_A(x)$ is a conjugacy invariant, so is $g_A(t)$. If $\alpha$ is a root of $\chi_A(x)$, then $\alpha+\alpha^{-1}$ is a root of $g_A(t)$.
Note that for  $\mu$ in $\C-\R$, $\mu +\mu^{-1}$ is a real number if and only if the norm of $\mu$ is $1$.

As in \secref{nature}, to detect the root of $g_A(t)$,  let
\begin{eqnarray*}
G'&=&(2a-c)-3.(\frac{-a}{3})(\frac{b-3}{3})+2(\frac{-a}{3})^3\\
&=&2a-c+\frac{1}{3}a(b-3)-\frac{2}{27}a^3\\
&=&\frac{1}{27}[27(a-c)+9ab-2a^3]=\frac{1}{27} G,
\end{eqnarray*}
\begin{eqnarray*}
H'&=&\frac{1}{3}(b-3)-\frac{1}{9}a^2\\
&=&\frac{1}{9}[3(b-3)-a^2]=\frac{1}{9}H.
\end{eqnarray*}
Define $\Delta'=G'^2 + 4H'^2$, and $\Delta=3^6 \Delta'=G^2 + 4H^3$.
Then $\Delta$ is the discriminant of the cubic equation \eqnref{1}. The multiplicity of a root of
$g_A(t)$ is determined by the resultant $R(g, g'')$ of $g_A(t)$ and its second derivative $g''_A(t)=6t-2a$. We observe that
$$R(g, g'')=-8[27(a-c)+9ab-2a^3]=-8G.$$
Now consider an element $A$ in $U(2,1;\F)$.

\medskip {\it Suppose $A$ is a hyperbolic.} Then $A$ is conjugate to an
element of the form (\ref{loxoe}) and  $g_A(t)$ has the following roots:
$$t_1=re^{{\bf i} \beta}+r^{-1}e^{{-\bf i} \beta}, t_2=r^{-1}e^{{\bf i} \beta}+re^{{-\bf i}, 
\beta}, t_3=2\cos\theta.$$

\noindent (1) If $\beta \neq 0$, then $g_A(t)$ has one real root, and two non-real complex conjugate
 roots $\lambda + \lambda^{-1}$ and $\bar{\lambda}+ \bar{\lambda^{-1}}$. Hence the norm of $\lambda$ must be different from $1$.
 Thus $A$ must be a regular hyperbolic. It follows from \secref{nature} that $\Delta>0$ and $G \neq 0$ in this case.

\medskip \noindent (2) Suppose $\beta=0$, or $\pi$, i.e.  $A$ acts a strictly hyperbolic, or a screw hyperbolic.  In this case $g_A(t)$ has a real root $r+r^{-1}$, or $-(r+r^{-1})$,  of multiplicity $2$. For $r>0$, note that $r+r^{-1}>2$.  Hence $g_A(t)$ has two and only two equal roots. Consequently we have $\Delta=0, \ G \neq 0$.
\subsubsection*{Suppose $A$ is elliptic}\label{elliptic} Then $A$ is conjugate to an element of the form (\ref{ellie}) and $g_A(t)$ has the following
roots:
$$t_1=2\cos\theta, t_2=2\cos\phi, t_3=2\cos\psi.$$

\noindent (3) If $A$ is regular elliptic, then all the roots of $g_A(t)$ are real and distinct from each-other. Hence it follows from \secref{nature} that $\Delta<0$.

\noindent (4) If $A$ is a complex elliptic or simple elliptic, then $g_A(t)$ has at least one repeated root, hence $\Delta=0$.
Note that $A$ is a simple elliptic if and only if $g_A(t)$ has a triple root, i.e.  $R(g, g'') = 0$.
Thus $A$ is a simple elliptic if and only if  $G = 0$. Hence if $A$ is complex elliptic, then $G \neq 0$. 

\subsubsection*{Suppose $A$ is parabolic}\label{parabolic} Then it is conjugate to an
element of the form (\ref{parae}). Hence all the roots of $g_A(t)$ are real. They are given by
$$t_1=2\cos\theta, t_2=2\cos\theta, t_3=2\cos\phi.$$
Since at least two of the roots are equal, hence  $\Delta=0$.

\noindent (5) If $A$ is screw parabolic, then $g_A(t)$ has a double root which is distinct from the third one, hence $G \neq 0$.

\noindent (6) Suppose $A$ has a triple root, then $A$ acts as either a unipotent, an ellipto-translation, or  an ellipto-parabolic.  Thus $g_A(t)$ also has a triple root, and hence  $G=0$. In this case, $|a|\leq 6, \ |b| \leq 15, \ |c| \leq 20$. If $A$ acts as a unipotent,  then $\theta=0=\phi$,  or $\theta=\pi=\phi$. Hence
$|a|=6,\ |b|=15,\ |c|=20$. Otherwise we have $|a|<6$, $|b|<15$,
$|c|<20$. Further note that the only non-parabolic class which has a triple root is the class of simple elliptics.  But for $g$ to be a simple elliptic, we must have $|a|<6, |b|<15, |c|<20$. The equality holds if and only if $\theta=0, \pi$, i.e. $g$ acts as the identity. Hence a non-trivial  isometry is unipotent if and only if $|a|=6,\ |b|=15,\ |c|=20$.

From the representative of the conjugacy class of $A$ it is clear that the elliptic and hyperbolic elements are semisimple.
If $A$ is simple elliptic, then the minimal polynomial $m_A(x)$ of $A_{\C}$ is linear. If $A$ is complex elliptic, then $m_A(x)$  is of the form \hbox{$(x-\lambda)(x-\mu)$}, $\lambda, \mu \in \s^1$. Hence the degree of the minimal polynomial is $2$.
 Parabolic isometries are not semisimple. We see that $A$ is an ellipto-translation or an ellipto-parabolic according as its
  minimal polynomial is of the form $(x-\lambda)^2$ or $(x-\lambda)^3$,  where $\lambda \in \s^1$. In particular,  $A$ is a
   vertical or a non-vertical translation according as the minimal polynomial is $(x-1)^2$ or $(x-1)^3$. If $A$ is a screw parabolic, then $m_A(x)$ is of the form $(x-\lambda)^2(x-\mu), \ \lambda, \mu \in \s^1$. 

Note from above that  $A$ acts as either screw hyperbolic,  strictly hyperbolic, screw parabolic, or  complex elliptic if and only if  $\Delta=0$, $G \neq 0$.  We further distinguish  these classes as follows. 
\begin{itemize}
\item[]{Let $\alpha_n=trace(A_{\C}^n)$.  
Thus  $\alpha_n$ is the sum of $n$-th powers of roots of $\chi_{A}(x)$.   In particular, $\alpha_1=a$.  For $A$ hyperbolic, $\alpha_n=2(r^n+r^{-n})+2 \cos^n \theta$, or, 
 $2. (-1)^n(r^n+r^{-n})+2 \cos^n \theta$. 
Thus 
$$|\alpha_n| \geq 2(r^n+ r^{-n})-2 |\cos^n \theta| \geq 2r^n-2.$$
Since $r>1$, we can choose $m$ such that $r^m>8$, i.e.  $|\alpha_m|>6$. Hence it follows that $\{\alpha_n\}$ is a divergent sequence.  Further $A$ acts as a strictly hyperbolic if and only if $\theta=0$ or $\pi$. In this case $g_A(t)$ must have a root $2$ or $-2$, and this implies from \eqnref{1} that $(2a+c)^2=4(1+b)^2$. For $A$ elliptic or parabolic, $\alpha_n=2(2\cos^n \theta+\cos^n \phi)$, and hence for all $n$,
 $|\alpha_n|<6$.   
}
\end{itemize}
\hspace{1 in} This completes the proof of \thmref{mainth}.

\section{The Centralizers and the $z$-classes: proof of \thmref{zcl}}\label{zcc}
Let $Z(x)$ denote the centralizer of $x$ in $U(2,1;\F)$, and
$Z_G(x)$ denote the centralizer of $x$ in the subgroup $G$ of
$U(2,1;\F)$. Let $T$ be an isometry of $\V^{2,1}$. An invariant
(right) subspace $\W$ of $\V^{2,1}$ is called
\emph{$T$-indecomposable} if
 it can not be expressed as a direct sum of two proper $T$-invariant (right) subspaces. Let $\oplus$ denote the orthogonal sum.

\subsection{Elliptic isometries} (1) Consider $E=E(\theta, \phi, \psi)$ with $\theta \neq \phi \neq \psi\neq \theta$. In this case $\V^{2,1}$ has a decomposition
into one dimensional $E$-invariant subspaces: $\V^{2,1}=l_1 \oplus
l_2 \oplus l_3$. Hence $Z(E)=Z(E|_{l_1}) \times Z(E|_{l_2}) \times
Z(E|_{l_3})$. This implies, $Z(E)=Z_{U(1;\F)}(e^{i\theta}) \times
Z_{U(1;\F)}(e^{i \phi}) \times Z_{U(1;\F)}(e^{i \psi})$. When
$\theta \neq 0$, the centralizer of $e^{i \theta}$ in $U(1;\F)$ is
the group $\s^1$ of all complex numbers of unit modulus. Hence
$$Z(E)=\left\{\begin{array}{ll} \s^1 \times \s^1 \times \s^1 \;\hbox { if }\theta, \phi, \psi \neq 0,\pi\\
\s^1 \times \s^1 \times \F^1 \hbox{ if one of }\theta, \phi, \psi=0 \hbox{ or }\pi\\
\F^1 \times \F^1 \times \s^1 \hbox{ if any two of }\theta, \phi,
\psi \hbox{ belongs to the set }\{0, \pi\}.
\end{array} \right.$$
where $\F^1$ is the multiplicative group consisting of elements of
$\F$ of unit modulus. When $\F=\C$, $\F^1=\s^1$, when $\F=\H$,
$\F^1=\s^3$.

\noindent (2) Let $\theta=\phi=\psi$. When $\F=\C$, $E=e^{i
\theta}I$ is an element in the center of $U(2, 1)$, hence
$Z(E)=U(2,1)$. When $\F=\H$,
$$Z(E)=\left\{\begin{array}{ll} U(2,1) \;\hbox { if }\theta \neq 0,\pi\\
Sp(2,1) \hbox{ otherwise.}
\end{array} \right.$$

\noindent (3) Suppose $\phi=\psi \neq \theta$. In this case,
$$Z(E)=Z_{U(1;\F)}(e^{i \theta}) \times Z_{U(2;\F)}(e^{i\phi}I_2).$$
When $\F=\C$, $Z(E)=\s^1 \times U(2)$.  When $\F=\H$,
$$Z(E)=\left\{\begin{array}{ll} \s^1 \times U(2) \;\hbox{ if } \phi \neq 0, \pi, \theta \neq 0,\\
\s^1 \times Sp(2) \;\hbox{ if } \theta \neq 0, \phi=0, \pi,\\
\s^3 \times U(2) \;\hbox{if }\theta=0, \pi, \phi \neq 0, \pi,\\
\s^3 \times Sp(2) \;\hbox{ if } \theta=0, \hbox{ resp. } \pi,
\phi=\pi, \hbox{ resp. }0
\end{array} \right.$$
Suppose $\phi=\theta \neq \psi$ (or $\psi=\theta \neq \phi$).  In
this case we have an orthogonal decomposition $\V=\W_2 \oplus \W_1$,
where $\dim \W_2=2$, $E|_{\W_2}=e^{i \theta} I_2$,
$E|_{\W_1}(v)=e^{i \psi}v$, where $v$ generates the one dimensional
subspace $\W_1$ and restriction of $\langle, \rangle$ on $\W_2$ has
signature $(1,1)$.

When $\F=\C$, $Z(E)=U(1,1) \times \s^1$.  When $\F=\H$,
$$Z(E)=\left\{\begin{array}{ll} U(1,1) \times \s^1 \ \hbox{ if }\theta \neq 0,\\
Sp(1,1) \times \s^1 \ \hbox{ if }\theta=0
\end{array} \right.$$
\noindent Given any elliptic element $T$, up to conjugacy, there are
the above choices for $Z(T)$. Hence there are four $z$-classes of
elliptic elements in $U(2,1)$,
 and eleven $z$-classes of elliptic elements in $Sp(2,1)$. In this count we have taken the identity map as an elliptic element.

\subsection{Hyperbolic isometries} Consider $L=L(\beta, \theta)$. In this case $L$ decomposes $\V^{2,1}$ orthogonally into $L$-invariant subspaces:
$\V^{2,1}=\V_2 \oplus \V_1$, where $\V_2$ is a $L$-indecomposable
subspace of $\V^{2,1}$ and restriction of $\langle, \rangle$ on
$\V_2$ has signature $(1,1)$. Hence $Z(L)=Z_{U(1,1;\F)}(L|_{\V_2})
\oplus Z_{U(1;\F)}(L|_{\V_1})$. When $\F=\H$,
$Z_{Sp(1)}(L|_{\V_1})=\s^1$ or $Sp(1)$ according to $\theta \neq 0$,
or $\theta=0$. For $\F=\C$, $Z_{U(1)}(L|_{\V_1})=\s^1$. Now, the
group $U(1,1;\F)$ is locally isomorphic to $SO_o(1,4)$, resp.
$SO_o(1,2)$ for $\F=\H$, resp. $\C$,  and $L$ acts as a hyperbolic
element. Thus the centralizers of $L|_{\V_2}$ is obtained from
\cite{gk}. Hence, given any hyperbolic element $T$, there is an
associated orthogonal decomposition $\V^{2,1}=\W_2 \oplus \W_1$ into
$T$-invariant subspaces as above, and we have
 $$Z(T)=Z_{U(1,1;\F)}(T|_{\W_2}) \times Z_{U(1;\F)}(T|_{\W_1}).$$

 When $\F=\C$, up to conjugacy, there is a unique choice for $Z_{U(1,1)}(T|_{\W_2})$, and $Z_{U(1)}(T|_{\W_1})$.
Hence there is a unique $z$-class of hyperbolic isometries of
$\h^2_{\C}$.

 It can be seen from the count in \cite{gk} that  there are two $z$-classes of hyperbolic elements in $SO_o(1,4)$. Hence, up to conjugacy, there are two choices
for $Z_{Sp(1,1)}(T|_{\W_2})$, and from the above we see that there are two choices for
$Z_{Sp(1)}(T|_{\W_1})$. Hence there are four $z$-classes of
hyperbolic isometries of $\h^2_{\H}$.

\subsection{Parabolic isometries}
Let $P$ be of the form (\ref{parae}).
 Let $A=\left(
                      \begin{array}{ccc}
                        a_1 & a_2 & a_3 \\
                        a_4 & a_5 & a_6 \\
                        a_7 & a_8 & a_9 \\
                      \end{array}
                    \right)\in \hat U(2,1;\F)$ be an element in $Z(P)$.
  It follows from $PA(\infty)=AP(\infty)=A(\infty)$ that
  $A\in G_{\infty}$, which implies that $a_2=a_3=a_8=0$, $\bar{a_1}a_5=1, |a_9|=1, \Re(\bar{a_1}a_4)=\frac{1}{2}|a_7|^2$ and $a_6=a_5\bar{a_7}a_9$.
  % and $a_6=a_5\bar{a_7}a_9$.
  By $AP=PA$, we get the following
  identities:
\begin{eqnarray}
\label{ineq1}&&  a_1e^{\bf i \theta}=e^{\bf i \theta}a_1,a_9e^{\bf i
\phi}=e^{\bf i
\phi}a_9,\\
\label{ineq2} && a_5f+a_6e^{\bf i \phi}=e^{\bf i \theta}a_6+fa_9, a_7e^{\bf i \theta}+a_9g=ga_1+e^{\bf i \phi}a_7,\\
\label{ineq3}&&a_4e^{\bf i \theta}+a_5d+a_6g=e^{\bf i
\theta}a_4+da_1+fa_7.
 \end{eqnarray}

\subsubsection{The unipotent isometries}  A unipotent isometry
is conjugate to
$$U=\begin{pmatrix}1 & 0 & 0\\s & 1 & \bar a\\a & 0 & 1\end{pmatrix},$$
where $Re(s)=\frac{1}{2}|a|^2$, and $a, s$ are complex numbers.

\noindent (1) $a=0$, i.e. $U$ is a vertical translation. Then $Z(U)$ is given by
  $$Z(U)=\bigg \{A=\left(
  \begin{array}{ccc}
    a_1 & 0 & 0 \\
    a_4 & a_1 & a_1\bar{a_7}a_9  \\
  a_7& 0 & a_9 \\
  \end{array}
\right) \ | \ A \in \hat U(2,1;\F), \  a_1\in \C \bigg \}.$$
 Hence there is a unique $z$-class of vertical translations in $Sp(2,1)$, resp. $U(2,1)$.

\noindent (2) $a \neq 0$, i.e. $U$ is a non-vertical translation. In this
case $\V^{2,1}$ is $U$-indecomposable. Let $S=U-I$, where $I$ is the
identity map. Let $V_U=(U-I)(\V^{2,1})$. Following  \cite{wall},
we call it the \emph{space} of $U$. Since $U$ is non-vertical translation,
kernel of $(U-I)$ is non-empty and is generated by a vector of norm
zero. Hence $V_U$ is of dimension $2$ over $\H$. The space $V_U$ can
be equipped with a non-degenerate Hermitian form $F_U$, cf. 
\cite[p.6]{wall},  given by
$$F_U(u, v) + \bar F_U(v, u)=\langle u, v \rangle. $$
We call it the \emph{form} of $U$.  The element $U$ is uniquely
determined by $V_U$ and $F_U$.
 The center $Z(U)$ is the isometry group of $(V_U, F_U)$.
Two non-vertical translations are in the same $z$-class if and
only if their forms are equivalent.

Now the nondegenerate Hermitian forms over $\H$ are determined up to
equivalence by their rank and signature,  cf. 
\cite[p.264]{lewis}. Since $\langle,\rangle$ has nonzero signature,
 hence the form of an isometry of $\H^{2,1}$ has nonzero signature. Now, there is a unique nondegenerate
Hermitian form of nonzero signature on a two dimensional vector
space over $\H$, viz. of type $(1,1)$. Hence, there is a unique
$z$-class of non-vertical translations in $Sp(2,1)$. By computation we see that
$$Z(U)=\bigg \{ A= \left(
  \begin{array}{ccc}
    a_1 & 0 & 0 \\
    a_4 & a_1 & a_6 \\
   a_7& 0 & aa_1a^{-1} \\
  \end{array}
\right) \ | \ A \in \hat U(2,1;\F), sa_1+\bar a a_7=a_1 s + a_6 a \bigg \}.$$

In particular, we note that the following elements commute with $U$:

(i) the vertical translations of the form: $\begin{pmatrix} 1 & 0 & 0 \\ s & 1 & 0 \\ 0 & 0 & 1 \end{pmatrix}$,  $s \in \F, \ Re(s)=0$,

(ii) the non-vertical translations of the form: $\begin{pmatrix} 1 & 0 & 0 \\ t & 1 & r\bar a\\ ra & 0 & 1 \end{pmatrix}, \ t \in \F, \ Re(t)=\frac{1}{2} |ra|^2$.

\subsubsection{The ellipto-translations and the ellipto-parabolics} Let $P$ be an ellipto-translation, resp.  ellipto-parabolic. Then by the Jordan decomposition, $P=P_s P_u$, where $P_s$ is a simple elliptic, and $P_u$ is a vertical, resp.  non-vertical translation. Hence $Z(P)=Z(P_s) \cap Z(P_u)$. We already know the descriptions of $Z(P_u)$, and that $Z(P_s)=U(2,1)$. Hence if $T_c$ is an ellipto-translation, then
$$Z(T_c)=\bigg \{ A= \left(
  \begin{array}{ccc}
    a_1 & 0 & 0 \\
    a_4 & a_1 & a_1\bar{a_7}a_9 \\
   a_7& 0 & a_9 \\
  \end{array}
\right) \ | \ A \in \hat U(2,1;\C)\bigg \}.$$
If $T_{\theta}$ is an ellipto-parabolic, then $Z(T_{\theta})$ is obtained similarly.

\subsubsection{The screw-parabolic elements}
We assume, without loss of generality, assume (for $\F=\H$)
 $$P=\begin{pmatrix} e^{i \theta} & 0 & 0 \\ d & e^{i \theta} & 0 \\  0 & 0 & e^{ i \phi} \end{pmatrix}, \ d\in \C-\{0\},
 \ \Re(e^{-i \theta}d)=0, \ 0 \leq \theta, \phi \leq \pi.$$
For $\F=\C$, $-\pi \leq \theta, \phi \leq \pi$. Then by the
equalities (\ref{ineq1})-(\ref{ineq3}), we have the following cases.

\noindent (i) Let $e^{i \theta},e^{ i \phi}\in \R$. We denote $P$ by
$P_o$ in this case.  It can be seen that
 $$Z(P_o)= \bigg\{A=\left(
  \begin{array}{ccc}
    a_1 & 0 & 0 \\
    a_4 & a_1 & 0 \\
   0& 0 & a_9 \\
  \end{array}
\right) \ | \ A \in \hat U(2,1;\F), \ a_1 \in \C \bigg\}.$$

\noindent (ii)  Let $e^{ i \theta} \in \R, \ e^{ i \phi}   \notin
\R$. Denote $P$ by $P_{o, \phi}$ in this case. Then
$$Z(P_{o,\phi})=\bigg\{A=\left(
  \begin{array}{ccc}
    a_1 & 0 & 0 \\
    a_4 & a_1 & 0   \\
    0& 0 & a_9 \\
  \end{array}
\right) \ | \ A \in \hat U(2,1;\F), \ a_1,a_9 \in \C \bigg\}.$$

\noindent (iii) Let $e^{\bf i \theta}\notin \R,e^{\bf i \phi}\in \R$.
Denote $P$ by $P_{\theta, o}$ in this case. Obviously, $a_1,a_5\in
\C$. By $a_4e^{\bf i \theta}+a_5d=da_1+e^{\bf i \theta}a_4$, we get
$a_4\in \C$, which in turn gives $a_1=a_5.$
 Therefore,
$$Z(P_{\theta, o})=\bigg\{A=\left(
  \begin{array}{ccc}
    a_1 & 0 & 0 \\
    a_4 & a_1 & 0   \\
    0& 0 & a_9 \\
  \end{array}
\right) \ | \ A \in \hat U(2,1;\F), \ a_1, a_4\in \C \bigg\}.$$

\noindent (iv)  Let $e^{\bf i \theta},e^{\bf i \phi}\notin \R$.
Denote $P$ by $P_{\theta, \phi}$ in this case. We see that
$$Z(P_{\theta, \phi})=\bigg \{ A= \left(\begin{array}{ccc}
    a_1 & 0 & 0 \\
    a_4 & a_1 & 0   \\
   0 & 0 & a_9 \\
  \end{array}
\right) \ | \ A \in \hat U(2,1;\C)  \bigg\}.$$

 Let $\F=\H$. From the description of the centralizers we see that there is a unique choice for the centralizers of ellipto-translations, a unique choice for the ellipto-parabolics, and there are four choices for the centralizer of a screw-parabolic.  Hence we have one $z$-class of ellipto-translations,
 one $z$-class of the ellipto-parabolics, and four $z$-classes of the screw-parabolic elements.

  When $\F=\C$, there is a unique $z$-class of
 screw-parabolic elements.

This completes the proof of \thmref{zcl}.

\subsection*{Acknowledgement}
 {\it  We appreciate Elisha Falbel, Ravi Kulkarni, John Parker and Ian Short for their comments on this work. John Parker gave us numerous valuable suggestions on some previous drafts of the paper.  Ian Short helped us to improve the English composition of the paper.  Special thanks to Parker and Short for their useful remarks. Finally it is a pleasure to thank the referee for many precious comments. }

\end{document}